\documentclass[review]{elsarticle}

\usepackage{hyperref}

\usepackage{latexsym} 
\usepackage{graphics}
\usepackage{epsfig}
\usepackage{graphicx}
\usepackage{tikz-qtree}
\usepackage{tikz}
\usepackage{amsmath}
\usepackage{dirtytalk} 
\usepackage{bbm} 

\usepackage{pgfplots}
\usepgfplotslibrary{groupplots}
\usepackage[utf8]{inputenc}
\usepackage[english]{babel}

\pgfplotsset{width=11cm,compat=newest}
\usepackage{algorithm}
\usepackage[noend]{algpseudocode}
\makeatletter
\def\BState{\State\hskip-\ALG@thistlm}
\makeatother
\algdef{SE}[DOWHILE]{Do}{doWhile}{\algorithmicdo}[1]{\algorithmicwhile\ #1}%

\newtheorem{proposition}{Proposition}

\journal{European Journal of Operational Research}

\biboptions{authoryear}

\begin{document}

\begin{frontmatter}

\title{Computing Optimal ($R,s,S$) Policy Parameters by a Hybrid of Branch-and-Bound and Stochastic Dynamic Programming}

\author{Andrea Visentin, Steven Prestwich}
\address{Insight Centre for Data Analytics, University College Cork, Ireland}
\ead{andrea.visentin@insight-centre.org,s.prestwich@cs.ucc.ie}

\author{Roberto Rossi}
\address{University of Edinburgh Business School, Edinburgh, UK}
\ead{roberto.rossi@ed.ac.uk}

\author{S. Armagan Tarim}
\address{Cork University Business School, University College Cork, Ireland}
\ead{armagan.tarim@ucc.ie}

\begin{abstract}
A well-know control policy in stochastic inventory control 
is the $(R,s,S)$ policy, in which inventory is raised to an
order-up-to-level \textit{S} at a review instant \textit{R} whenever it falls 
below reorder-level \textit{s}. To date, little or no
work has been devoted to developing approaches for computing $(R,s,S)$ policy parameters.  In this work, we introduce a hybrid approach that exploits tree search to
compute optimal replenishment cycles, and stochastic dynamic
programming to compute $(s,S)$ levels for a given cycle.  Up to 99.8\%
of the search tree is pruned by a branch-and-bound technique with
bounds generated by dynamic programming.  A numerical study shows that
the method can solve instances of realistic size in a reasonable time.
\end{abstract}


\end{frontmatter}


\section{Introduction}

Single-item, single-stocking location, stochastic inventory systems
have long been investigated under various operational assumptions, and
the associated literature is large.  Scarf's seminal paper, \cite{scarf1959optimality}, addressed this problem over a finite
planning horizon comprising discrete time periods, non-stationary
stochastic demands, a fixed ordering cost, and linear holding and
shortage costs.  Scarf proved that the $(s,S)$ policy (more precisely
the $(s_t,S_t)$ policy) is cost-optimal.  In this policy the decision-maker checks the current inventory position at review epochs (the
start of each time period) and if the inventory position is at or
below $s_t$ an order is placed to raise it to $S_t$.  For a planning
horizon of $T$ periods the optimal $(s_t,S_t)$ policy requires $2T$
policy parameters, computed in a here-and-now fashion at the start of
the planning horizon.  Actual replenishment timings and associated
order quantities are instead determined in a wait-and-see manner. 

In this paper, we address a more general form of the inventory control problem described by Scarf. 
According to \cite{silver1981operations} the $(R,s,S)$ policy is one
of the most commonly adopted inventory control strategies (also
called $(T,s,S)$ or $(s,S,T)$ in the literature, \cite{babai2010empirical, lagodimos2012computing}).  In an
$(R_t, s_t, S_t)$ system the inventory level is checked only at review epochs
$R_t$, which are policy parameters that are fixed at the start of the
planning horizon. After a review, an order may be placed to raise the inventory level up to $S_t$ if it is at or below $s_t$. 

Two special cases of the $(R,s,S)$ policy naturally arise.  Firstly,
it reduces to the $(s,S)$ case if there is no explicit cost involved
in carrying out inventory reviews. Inventory review (also known
as stock-taking) is costly in practice, so we consider the case in
which a fixed {\it system control\/} cost \cite{silver1981operations}
is incurred when the inventory is reviewed, e.g. \cite{fathoni2019development,christou2020fast}. The $(R,s,S)$ policy relaxes the cost accounting assumption that the fixed cost of replenishment covers both review and delivery costs, and separates the fixed cost of conducting a review from the fixed ordering cost. One practical implication of this relaxation is that the order cancellation and relevant costs can be explicitly incorporated into inventory planning.

Secondly, the $(R,s,S)$ reduces to the
$(R,S)$ policy (the {\it replenishment cycle policy\/}) if reorder
levels $s_t$ are equal to the order-up-to-levels $S_t$.  In a replenishment cycle policy,
the replenishment periods are fixed at the beginning of the planning
horizon and the replenishment orders are placed only in these periods
after period demands so far have been observed.

Although $(R,s,S)$ is one of the most general and frequently used
inventory policies, as pointed out by \cite{silver1981operations} {\it
  the determination of the exact best values of the three parameters
  is extremely difficult\/}.  To the best of our knowledge no approach
to computing them has been presented in the literature.  We fill this
important gap in the literature by making the following contributions:
\begin{itemize}
\item
we introduce an efficient hybrid of branch-and-bound and stochastic dynamic program (SDP) to
compute optimal policy parameters;
\item
we improve the branch-and-bound by using tighter bounds computed
through a separate dynamic programming (DP) method;
\item
we show empirically that the new algorithm performs significantly
better than a baseline method and that it is able to solve problems of realistic size in a reasonable time.
\end{itemize}
The paper is structured as follows.  Section
\ref{sec:literature_review} surveys related literature.  Section
\ref{sec:prob_description} provides a problem description.  Section
\ref{sec:sdp_ss} introduces a simple SDP formulation.  Section
\ref{sec:method} introduces a branch-and-bound strategy.  Section
\ref{sec:experimental} carries out a comprehensive numerical study.
Finally, Section \ref{sec:conclusion} concludes the paper.

\section{Literature review} \label{sec:literature_review}
The problem of computing policy parameters for an inventory control system under stochastic demand has received a great deal of attention. In this section, we survey the relevant literature on the classic stochastic inventory control problem. We then survey different versions of the problem. Finally, we survey $(R,s,S)$ real-world applications. 

An important class of these problems is the single-item single-location non-stationary stochastic lot-sizing under linear holding costs, penalty costs and both linear and fixed ordering costs. Different policies can be used to determine the size and timing of the orders on such setting.  
In his seminal work, Scarf characterises the structure of the optimal policy for such a problem. The framework proposed by \cite{bookbinder1988strategies} divides the policies into three classes: static uncertainty, dynamic uncertainty and static-dynamic uncertainty. These classes differ in the moment at which the decisions are taken. Since then, numerous research works tackled the computation of policy parameters under demand uncertainty, mainly focusing on the $(s,S)$ and the $(R,S)$ policies that have a flexible order quantity. According to the strategies categorisation presented in \cite{powell2019unified}, these works can be divided into two types: \textit{deterministic/special structure solutions} or \textit{sample models}. The first category, that comprises this study, includes a wide variety of approaches, based on: dynamic programming \cite{scarf1959optimality,rossi2011state,ozen2012static}, mixed-integer linear programming \cite{tarim2004stochastic,xiang2018computing,tunc2018extended}, approximations \cite{gutierrez2017simple}, and constraint programming \cite{rossi2012constraint}. The sample models category includes two-stage stochastic programming, that has been applied to inventory policy computation in \cite{fattahi2015investigating,cunha2017two,dos2019enhanced}. 
Different comparison studies have been conducted recently to benchmark different aspects of these policies: \cite{kilic2011investigation} extends a measure of planning instability for the non-stationary stochastic lot-sizing; \cite{dural2019benefit} compares different policies performances in the receding horizon, \cite{sani1997selecting} and \cite{babai2010empirical} are comparative studies on the performances of $(s,S)$ heuristics. 

Modifications on the original inventory model have been proposed to allow a closer representation of real-world problems. \cite{dillon2017two} proposes an $(R,S)$ policy solution to manage the blood supply chain, that includes perishable products. \cite{alvarez2020inventory}'s model considers both quantity and quality decay of the inventory product; the quality can be improved by mixing it with a higher quality product. A set of heuristics for the lot-sizing problem with remanufacturing of returned products is presented in \cite{kilic2019heuristics}. All-units discount $(s,S)$ policy has been analysed in \cite{wang2019procurement}. Uncertainty can involve other aspects of the inventory system; for example, \cite{bashyam1998optimization,rossi2010computing} considers a stochastic lead time. Different supply chain configuration can be considered; for example, a two-echelon inventory system \cite{schneider1991empirical,schneider1995power}. \cite{ma2019stochastic} provides an updated review on stochastic inventory control algorithms, while \cite{bushuev2015review} presents a broader picture of the state-of-the-art in lot sizing.

The $(R, s, S)$ policy parameters computation has been tackled in the literature under the stationary, continuous time setting. With this configuration, only three parameters have to be optimised since the demand does not change over time. This problem has been solved to optimality by \cite{lagodimos2012computing}. In \cite{lagodimos2012optimal,christou2020fast} a batch version of the policy is considered.

None of the surveyed methods can be easily adapted to compute the $(R,s,S)$ policy parameters under the finite horizon and discrete time setting: since it has three sets of decision variables making the previous models not applicable. 
While other policies can be used for the same problem, the $(R,s,S)$ is a more generic model and has better cost performances, being the $(s,S)$ and the $(R,S)$ specific case of an $(R,s,S)$ policy. The introduction of the review cost makes no difference in the $(s,S)$ and $(R,S)$ policy computation; in the $(s,S)$ the cost is charged in every period, while in the $(R,S)$ every review coincides with an order. A static policy would have poor performance as well because it can not react to the demand realisations \cite{dural2019benefit}.


The $(R,s,S)$ policy is widely used by practitioners, usually not
independently but as a component of complex supply chains, and here we
survey some recent models.  Due to the
complexity of the determination of its parameters, in the surveyed
papers, the value of $R$ is considered to be constant across the time
horizon. \cite{bijvank2012inventory} describe an
inventory control system for point-of-use location.  They compare the
performance of $(R,s,Q)$ policies (with fixed order quantities)
against $(R,s,S)$ under stationary stochastic demand.  Because of
stationarity, the policy parameters were constant throughout the
horizon.  \cite{ahmadi2018optimal,monthatipkul2008inventory} tackle a
capacitated two-echelon inventory system with one warehouse and
multiple retailers.  They use a heuristic based on
\cite{schneider1995power} for the $(R,s_t,S_t)$ policy.
\cite{cabrera2013solving} consider a similar two-level supply chain in
which a single plant serves a set of warehouses, which in turn serve a
set of end customers or retailers.  The warehouses model is based on
$(R,s,S)$ and they develop a heuristic to solve an inventory location
model with this configuration.  The same problem has been tackled by
\cite{araya2018lagrangian} using Lagrangian relaxation and the
subgradient method.  \cite{bijvank2012lost} analysed lost-sales
inventory control policies with service level.  They define an optimal
policy starting from the $(s,S)$ SDP introduced by
\cite{scarf1959optimality}.  They present a value-iteration algorithm
to find the $(R,s,S)$ parameters that minimise the inventory cost
subjected to service constraints.  As the parameters are fixed, their
solution is unsuitable for a non-stationary setting.

The analysis of the state-of-the-art confirms the novelty of the solution, and the practitioners' interest in the $(R,s,S)$ policy usage in a stochastic environment.

\section{Problem description} \label{sec:prob_description}

We consider the single-item, single-stocking location, stochastic
inventory control problem over a $T$-period planning horizon.  Without
loss of generality, we assume that orders are placed at the start of
each period and that the lead time is zero, as is common in the
literature
\cite{scarf1959optimality,bollapragada1999simple,tarim2004stochastic}.
An inventory control policy defines the timing and quantities of
orders over the planning horizon.  We define a review moment, or
review period, as a period in which the level of the inventory is
assessed and an order can be placed.  A replenishment cycle is
represented by the interval between two review moments.
We denote by $Q_t$ the quantity of the order placed in period $t$, and
an inventory review cost by $W$.  Ordering costs are represented by a
fixed value $K$ and a linear cost, but we shall assume without loss of
generality that the linear cost is zero.  The extension of our
solution to the case of a non-zero production/purchasing cost is
straightforward, as this cost can be reduced to a function of the
expected closing inventory level at the final period
\cite{tarim2004stochastic}.  At the end of each period, a linear
holding cost $h$ is charged for every unit carried from one period to
the next.

Demands $d_t$ in each period $t$ are independent random variables with known probability distributions.  Backlogging of excess demand is
assumed, so if the demand in a period exceeds on-hand inventory the
rest of the demand is carried to the next period; a linear penalty
cost $b$ is incurred on any unit of back-ordered demand at the end of
a period.

Under the non-stationarity assumption the $(R,s,S)$ policy takes the
form $(R_t,s_t,S_t)$ where $R_t$ denotes the length of the $t^{th}$
replenishment cycle, while parameters $s_t$ and $S_t$ denote the
reorder-level and order-up-to-level associated with the $t^{th}$
inventory review.  We consider the problem of computing the $(R,s,S)$
policy parameters that minimize the expected total cost over the
planning horizon.

\section{Stochastic dynamic programming formulation} \label{sec:sdp_ss}

In this section, we provide a simple technique to compute the optimal
$(R,s,S)$ policy parameters.  It can be considered the
state-of-the-art in computing such parameters in the presence of
stochastic non-stationary demand.  Moreover, it constitutes the basis
of the branch-and-bound technique introduced later.

We represent the replenishment moments by binary variables $\gamma_t$
($t = {1, \dots, T}$) which take value 1 if a review is placed in
period $t$ and 0 otherwise.  We assume $Q_t = 0$ if $\gamma_t = 0$ so
no order will be placed outside a review moment.  The optimal
$(R,s,S)$ policy for our problem is represented by the parameters
$\gamma_t,s_t,S_t$ that minimize the expected total cost.

Consider an arbitrary review cycle plan with $\gamma_t$ as a
parameter, not a decision variable.
We denote the closing inventory level for each period by $I_t$, and
the given initial inventory level by $I_0$.  We assume that the orders
are placed at the beginning of each time period and delivered
instantaneously.  The problem can be formulated and solved to
optimality as an SDP \cite{bellman1966dynamic}.
The expected immediate cost combining ordering, review, holding and penalty costs,
given action $Q_t$:
\begin{align}
\label{eq:immediate_cost}
f_t(I_{t-1}, Q_t) &=& \gamma_t W +  K \mathbbm{1}\{ Q_t>0 \} +
E[h \max(I_{t-1} + Q_t -d_t, 0)  +\nonumber\\
&& b \max(d_t-I_{t-1}- Q_t, 0)]
\end{align}

Let $C_t(I_{t-1})$ represent the expected total cost of an optimal policy over periods $t, \dots, \ n$ and $\mathbbm{1}$ is the indicator function. These variables are the states of the DP formulation. We model the problem with the functional equation:
\begin{equation}
\label{eq:functional_eq}
C_t(I_{t-1}) =  \underset{0 \leq Q_t \leq M \gamma_t}{\min} ( f_t(I_{t-1}, Q_t) +
E[ C_{t+1}(I_{t-1}+ Q_t - d_t) ]  ) 
\end{equation}
where $M$ is a sufficiently large number.
The boundary condition is:
\begin{equation}
C_{T+1}(I_{T}) = 0
\end{equation}

$C_{1}(I_0)$, where $I_0$ is the initial inventory level, contains the
expected cost for the optimal $(s,S)$ policy associated with the
\textbf{$\gamma$} assignment.  To reduce the computational time we can
exploit the property of K-convexity \cite{scarf1959optimality} when
solving the SDP.

Let $\hat{C}_1(I_0)$ represent the expected total cost of the optimal
$(R,s,S)$ policy, given the initial inventory level $I_0$ at period 1.
We can define it as:
\begin{equation}
\label{eq:baseline}
\hat{C}_1(I_{0}) = \underset{\gamma_1, \dots, \gamma_T}{\min} ( C_1(I_0) ) 
\end{equation}
Evaluating the optimal $(s,S)$ policy for all possible assignments of
$\gamma_1, \dots, \gamma_T$ yields the optimal $(R,s,S)$ policy. The model works with every possible demand distribution, as long as it is finite and discretisable. This
is our baseline method on which we aim to improve.

\subsection{Unit cost} \label{sec:unit_cost}
The algorithm can be extended to model the per unit ordering cost. There are two options, reducing it to a function of the expected closing inventory, e.g. \cite{tarim2004stochastic}; or including it in the immediate cost function. 

Let $v$ be the per unit ordering/production cost, Equation \ref{eq:immediate_cost} is replaced by:
\begin{align}
f_t(I_{t-1}, Q_t) =& \gamma_t W +  K \mathbbm{1}\{ Q_t>0 \} +  v  Q_t +
\nonumber\\
& E[h \max(I_{t-1} + Q_t -d_t, 0)  + b \max(d_t-I_{t-1}- Q_t, 0)]
\end{align}

\subsection{Lost sales} \label{sec:lost_sales}
Complete backlogging of the demand is a limiting assumption in many real-world settings. Studies analysing customer behaviour show that in case of a stock out, only a minority delay the purchase \cite{verhoef2006out}. According to \cite{bijvank2012inventory}, the lost sales configuration is underrepresented in the lot-sizing literature, even if it is more appropriate to model customers' behaviour. Approximating a lost sales model with a backlog model results in a non-negligible increase in costs \cite{zipkin2008old}.

The SDP formulation can be extended to model lost sales. We considered the partial backorder configuration presented in \cite{dos2019enhanced}. They define as $\beta$ ($\beta \in [0,1]$) the fraction of the unmet demand that is carried on to the next period and the reminder is lost. This parameter gives the flexibility to model both backlog ($\beta=1$), lost sales ($\beta= 0$) or a combination of the two. The functional equation \ref{eq:functional_eq} becomes:
\begin{equation}
\label{eq:functional_eq_lost}
C_t(I_{t-1}) =  \underset{0 \leq Q_t \leq M \gamma_t}{\min} ( f_t(I_{t-1}, Q_t) +
E[ C_{t+1}(\max(I_{t-1}+ Q_t - d_t, \beta (I_{t-1}+ Q_t - d_t)) ]  ) 
\end{equation}

\subsection{Example} \label{example}

We use a simple example to illustrate the application of our method,
with a 3-period planning horizon.  We assume an initial inventory
level of zero and a Poisson distributed demand for each period with
averages $\overline{d} = [20,30,40]$.
We consider an ordering cost value $K=30$, a review cost $W=10$, and
holding and penalty costs of $h=1$ and $b=10$ per unit per period
respectively.

The algorithm must choose replenishment moments $\gamma = \langle
\gamma_1, \gamma_2, \gamma_3 \rangle$ that minimize the expected cost
of the policy.  Table \ref{tab:toy_baseline} shows the expected cost
of each $(s,S)$ policy computed with different review periods.  The
optimal solution is $\gamma = \langle 1,0,1 \rangle$ with expected
cost 142.7.  However, exhaustive search becomes impractical as the
planning horizon grows so in Section \ref{sec:method} we develop a
more efficient method.

\begin{table*}
\begin{center}
\begin{tabular}{|c|c|c||c|}
\hline
$\gamma_1$ & $\gamma_2$ & $\gamma_3$ & Expected cost \\ \hline
\hline
0 & 0 & 0 & 1600.0 \\ \hline
0 & 0 & 1 & 751.8 \\ \hline
0 & 1 & 0 & 304.7 \\ \hline
0 & 1 & 1 & 302.0 \\ \hline
1 & 0 & 0 & 185.0 \\ \hline
1 & 0 & 1 & \textbf{142.7} \\ \hline
1 & 1 & 0 & 153.1 \\ \hline
1 & 1 & 1 & 150.4 \\ \hline
\end{tabular}
\end{center}
\caption{Optimal expected cost for the 3-period example.}
\label{tab:toy_baseline}
\end{table*}

\section{A hybrid of branch-and-bound and SDP} \label{sec:method}

In this section, we present a hybrid technique that combines SDP and
branch-and-bound.  The algorithm obtains optimal $(R,s,S)$ policies
associated with specific review plans at leaf nodes.  The search tree
(defined in Section \ref{sec:searchtree}) is explored by depth-first
search (DFS).  The subproblems associated with the nodes are defined
in Section \ref{sec:subproblems}.  Section \ref{sec:pruning}
introduces the pruning condition and lower bound computed with DP.
Finally, Section \ref{sec:nodes} presents the node resolution process.

\subsection{Search tree} \label{sec:searchtree}

\begin{figure}
\tikzset{every tree node/.style={minimum width=2em,draw,circle},
         blank/.style={draw=none},
         edge from parent/.style=
         {draw,edge from parent path={(\tikzparentnode) -- (\tikzchildnode)}},
         level distance=2cm,
        level 3/.style={sibling distance = 1cm}
        }
\centering

\begin{tikzpicture}
\Tree
[.  1 
    \edge node[auto=right] {$\gamma_3 = 1$};  
    [.  2 
        \edge node[auto=right] {$\gamma_2 = 1$};  
        [.  3
            \edge node[auto=right] {$\gamma_1 = 1$};  
            [.4 ]
            \edge node[auto=left] {$\gamma_1 = 0$};
            [.4 ]
        ]
        \edge node[auto=left] {$\gamma_2 = 0$};  
        [.3 
            [.4 ]
            [.4 ]
        ]
    ]
    \edge node[auto=left] {$\gamma_3 = 0$};
    [.2  
        \edge node[auto=right] {$\gamma_2 = 1$};  
        [.3 
            [.4 ]
            [.4 ]
        ]
        \edge node[auto=left] {$\gamma_2 = 0$};  
        [.3 
            [.4 ]
            [.4 ]
        ]
    ]
]
\end{tikzpicture}
\caption{Search tree for a 3-period instance: nodes contains level
  numbers.}
\label{fig:binary_tree}
\end{figure}
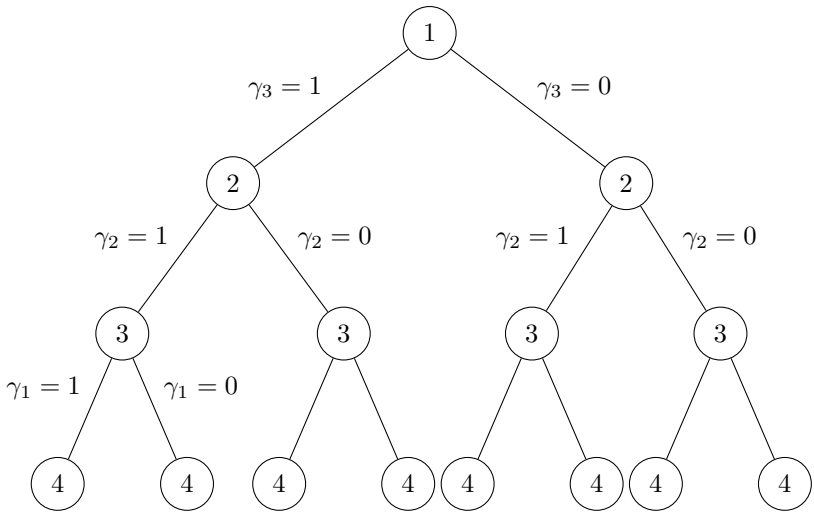

The branch-and-bound goal is to find the review plan with the minimum
expected cost.  During branching of $\gamma_t$, the value is fixed to $1$ or
$0$.  The search tree has $T+1$ levels, and the branching at its root
fixes the value of $\gamma_T$.  At level $\ell$ branching involves the
variable $\gamma_{T- \ell+1}$.  The path from the root to a node at
level $\ell$ represents a fixed assignment of the suffix $\langle
\gamma_{T- \ell +2}, \dots, \gamma_{T} \rangle$.  A leaf node
represents a complete assignment of the $\gamma$ values.  Figure
\ref{fig:binary_tree} shows the search tree of a 3-period problem, as
in the example presented in the previous section.

\subsection{Subproblems} \label{sec:subproblems}

Given the period $t$ and the partial assignment of a suffix of the
review moments $\langle \gamma_t, \dots, \gamma_T \rangle$, the
problem at a node is to find the $\langle \gamma_1, \dots \gamma_{t-1}
\rangle$ that minimizes the expected cost of the optimal policy.  We
denote this problem as BnB-SDP($t$,$\langle \gamma_t, \dots, \gamma_T
\rangle$).  For each subproblem using Equation \ref{eq:functional_eq}
we can compute the expected cost of the optimal policy starting at
period $t$ with inventory level $i$.  This is possible because all
review moments after period $t$ are fixed, and because of the SDP
stage structure presented in Section \ref{sec:sdp_ss}.

\subsection{Bounds and pruning} \label{sec:pruning}

If all the solutions in the subtree rooted in a node are suboptimal
then we can prune that node without compromising
optimality.
\begin{proposition}
\label{prop:property_monotonical}
Given a fixed assignment of \textbf{$\gamma$}:
\begin{equation}
\label{eq:property_monotonical}
 \underset{I}{\min}(C_t(I)) \geq \underset{I}{\min}(C_{t-1}(I))
\end{equation}
\end{proposition}
From the functional equation (\ref{eq:functional_eq}) it is clear that
$C_t$ is equal to the expected value of $C_{t+1}$ plus some
non-negative costs, so the minimum cost in each stage increases
monotonically with tree depth.

During tree search $\bar{C}$ records the expected cost of the best
plan computed so far, that is the minimum $C_1(I_0)$ among all leaves
already computed.  This is used as an upper bound for the expected
cost of the optimal plan as follows.  Considering the subproblem
BnB-SDP($t$,$[\gamma_t, \dots, \gamma_T]$) with the associated
$C_t(i)$ expected costs:
\begin{proposition}
If
\begin{equation}
\label{eq:pruning_condition}
\underset{i}{\min}(C_t(i)) \geq \bar{C}
\end{equation}
then because of the monotonicity of the cost function:
(\ref{eq:property_monotonical}):
\begin{equation}
\underset{i}{\min}(C_1(i)) \geq \bar{C}
\end{equation}
Finally, since the expected cost associated with a plan ($C_1(I_0)$)
is part of $C_1$:
\begin{equation}
C_1(I_0) \geq \bar{C}
\end{equation}
\end{proposition}
Hence if (\ref{eq:pruning_condition}) is true the subproblem
BnB-SDP($t$,$[\gamma_t, \dots, \gamma_T]$) is not part of an optimal
solution and the search tree can be pruned.

However, this pruning condition makes no assumption on the costs faced
on periods $ 1, \dots , t-1$, and a lower bound on the costs in those
periods leads to more effective pruning.  Let $MC_t(I_t)$ represent a
lower bound on the cost faced in periods $1, \dots , t$ with a closing
inventory of $I_t$ in period $t$.  The pruning condition
(\ref{eq:pruning_condition}) can be refined to:
\begin{equation}
\label{eq:pruning_condition_updated}
\underset{I_t}{\min}(C_t(I_{t-1}) + MC_{t-1}(I_{t-1}) ) \geq \bar{C}
\end{equation}
Having a bound independent from the review periods allows us to compute
it only once before the branch-and-bound algorithm.

The bounds can be computed by a DP with stages and states equivalent
to the SDP presented in Section \ref{sec:sdp_ss} and functional
equation:
\begin{equation}
MC_t(I_{t}) = \min
\left\{\begin{matrix}
  f_t(I_{t}, 1) + \underset{j < I_t}{\min} (MC_{t-1}(j))\\ 
  f_t(I_{t}, 0) + \underset{j \geq I_t}{\min} (MC_{t-1}(j))
\end{matrix}\right.
\label{eq:funceqdp}
\end{equation}
where, as defined in Section \ref{sec:sdp_ss}, $I_t$ is the current
inventory level, and $ f_t(I_{t}, Q_t)$ is the
ordering-holding-penalty cost.  In the first case, an order has been
placed in period $t$ so the inventory level in the previous period was
less than or equal to the current level.  In the second case, an order
has not been placed so the previous inventory level was greater than
or equal to the current level.  The boundary condition is:
\begin{equation}
MC_1(I_{1}) =
\left\{\begin{matrix}
  W + K  + f_1(I_{1})  & \text{if } I_1 > I_0\\ 
  f_1(I_{1}) & \text{if } I_1 \leq I_0
\end{matrix}\right.
\end{equation}
where $I_0$ is the initial inventory. Considering finite demand, the DP has an amount of states equal to the number of periods multiplied by the maximum inventory level. Each state requires a single  computation of Equation \ref{eq:funceqdp}, that is pseudo-polynomial in relation to the maximum inventory. The overall complexity of a DP is the number of states multiplied by the complexity required to solve one of them, so the overall complexity is pseudo-polynomial.

\subsection{Node computation} \label{sec:nodes}

Algorithm \ref{alg:bnb} summarises the branch-and-bound procedure.  In
line 1, the SDP stage $t$ is solved.  In line 7, the pruning condition
is evaluated: if a pruning occurs the branching phase is skipped.  In
lines 8 and 9 DFS recursively continues.  Lines 3--6 relate to leaf
nodes: if the policy represented by the leaf is better than the best
found so far, the value of $\bar{C}$ is updated.  The algorithm starts
by invoking BnB-SDP($T+1$,$\emptyset$), and at the end, the expected
cost of the optimal policy is given by $\bar{C}$.

The algorithm as always shown branches by assigning first $\gamma_t =
0$, but its performance can be improved by randomisation.  If, during
each branching phase, we randomly order lines 8--9 we obtain a better
solution earlier, leading to a stronger pruning of the search tree.
We evaluate the effect of this randomisation in Section
\ref{sec:experimental}.

\begin{algorithm}
\caption{BnB-SDP($t$,$[\gamma_t, \dots, \gamma_T]$)}\label{alg:bnb}
\emph{Data}: the current upper bound $\bar{C}$, the $C_{t+1}(i)$
computed at the parent node, and the bounds $MC(i)$.
\begin{algorithmic}[1]
\State \text{Compute $C_t$ using Equation \ref{eq:functional_eq}}
\If {$t = 1$}
    \If {$C_1(I_0) < \bar{C}$}
        \State $\bar{C} \gets C_1(I_0)$
        \State Save $[\gamma_1, \dots, \gamma_T]$ as incumbent review plan
    \EndIf
\Else
    \If {$\min(C_t(i) + MC_{t-1}(i)) \geq \bar{C} $}
        \Return
    \EndIf
    \State BnB-SDP($t-1$,$[0, \gamma_t, \dots, \gamma_T]$)
    \State BnB-SDP($t-1$,$[1, \gamma_t, \dots, \gamma_T]$)
\EndIf
\end{algorithmic}
\end{algorithm}

\subsection{Guided tree search} \label{sec:guided}
The random descent can get stuck in inferior branches of the search tree. It takes a considerable time to get a reasonable review plan, and a good cost upper bound in these cases. Computing a near-optimal review plan and using it to guide the search leads to the immediate computation of a policy with a low expected cost. This tighter bound increases the number of nodes proved not-optimal by the pruning condition.

A reasonable review plan can be computed using the $(R_t,S_t)$ policy. As mentioned in the introduction, this policy places an order at each review moment. The replenishment cycles ($R_t$) can be used as a review plan, while the order-up-to-levels $S_t$ can be ignored. During the first descend of the branch-and-bound search tree, the delta values are selected following this review plan; thus, the first leaf to be computed is the one that has $R_t$ as review moments. This leaf represents the optimal $(R_t, s_t, S_t)$ policy for that review plan, and it should have a low expected cost.
After computing the first leaf of the tree, the search proceeds in the replenishment plan's neighbourhood using a randomised approach.

The experimental section shows the improvement in pruning efficacy and computational time. 
Good computational performance and implementation simplicity make the MILP formulation presented in \cite{rossi2015piecewise} a good solution to compute the $(R_t,S_t)$ policy; this formulation is used in the experimental section.

\subsection{Example} \label{sec:example_pruning}

\begin{figure}
\tikzset{every internal node/.style={minimum width=2em,draw,circle},
        every leaf node/.style={minimum width=2em,draw,circle,dashed},
         blank/.style={draw=none},
         edge from parent/.style=
         {draw,edge from parent path={(\tikzparentnode) -- (\tikzchildnode)}},
         level distance=2cm,
        level 3/.style={sibling distance = 1cm}
        }
\centering

\begin{tikzpicture}
\Tree
[.  0
    \edge node[auto=right] {$\gamma_3 = 1$};  
    [.107
        \edge node[auto=right] {$\gamma_{2} = 1$};  
        [.142
            \edge node[auto=right] {$\gamma_{1} = 1$};  
            [.\textbf{150} ]
            \edge node[auto=left] {$\gamma_{1} = 0$};
            [.302 ]
        ]
        \edge node[auto=left] {$\gamma_{2} = 0$};  
        [.138
            \edge node[auto=right] {$\gamma_{1} = 1$};  
            [.\textbf{143} ]
            \edge node[auto=left] {$\gamma_{1} = 0$};  
            [.752 ]
        ]
    ]
    \edge node[auto=left] {$\gamma_3 = 0$};
    [.109
        \edge node[auto=right] {$\gamma_{2} = 1$};  
        [.\underline{144}
            \edge[dashed];
            [.* ]
            \edge[dashed];
            [.* ]
        ]
        \edge node[auto=left] {$\gamma_{2} = 0$};  
        [.\underline{181}
            \edge[dashed];
            [.* ]
            \edge[dashed];
            [.* ]
        ]
    ]
]
\end{tikzpicture}
\caption{Branch-and-bound technique applied to the toy problem.}
\label{fig:binary_tree_toy_pruning}
\end{figure}
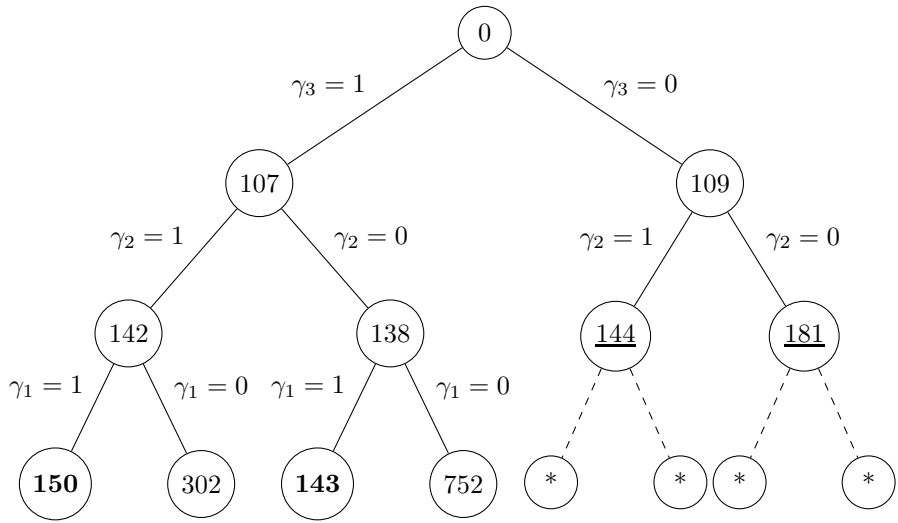

The search tree with the DP bounds for the example of Section
\ref{example} is represented in Figure
\ref{fig:binary_tree_toy_pruning}.  Each internal node contains the
value of the pruning condition with the DP bounds
(\ref{eq:pruning_condition_updated}).  An internal node is underlined
if the pruning occurs in that node.  Each leaf is in bold if it
contains an improvement compared to the previous best solution
$\bar{C}$.  Pruned nodes are indicated by an asterisk (*).

We define \textit{pruning percentage} as the percentage of nodes that
are proved to be suboptimal by the pruning condition during tree
search.  In this example, the number of computed nodes is 10 and 4
nodes have been pruned, so the pruning percentage is $4/14=28.57\%$.

\section{Computational study} \label{sec:experimental}

In this section, we evaluate the new methods, including an
assessment of the effects of branching randomisation and problem
parameters (costs) empirically.  We conduct two sets of experiments as follows.
In Section \ref{sec:scalability}, we analyse the scalability of the new
approaches by increasing the number of periods until no method is able
to solve the problem within a 1-hour time limit consistently.  In
Section \ref{sec:type_analysis} we fix the planning horizon to 10 and 20 periods and vary the cost parameters.  
For the experiments, we use
three $(R,s,S)$ policy solvers:
\begin{itemize}
\item
\textbf{SDP}, the SDP technique described in Section \ref{sec:sdp_ss}
which we consider the current state-of-the-art.
\item
\textbf{BnB}, the branch-and-bound solution introduced in Section
\ref{sec:method}.
\item
\textbf{BnB-Rand}, branch-and-bound with randomised branching.
\item
\textbf{BnB-Guided}, branch-and-bound with a guided tree search, Section \ref{sec:guided}.
\end{itemize}
We compare these in terms of computational time, pruning percentage
and average number of review periods (but not expected costs because
the solutions are optimal in each case).  All experiments are executed
on an Intel(R) Xeon E5620 Processor (2.40GHz) with 32 Gb RAM.
For the sake of reproducibility, we made the code available\footnote{\url{https://github.com/andvise/RsS-EJOR}}.


We base our numerical studies on the set of instances originally
proposed by \cite{berry1972lot} and widely used in the literature
\cite{dural2016comparison,rossi2008global,xiang2018computing}.  A Poisson variable represents the demand in each period.

\subsection{Scalability} \label{sec:scalability}
This experiment aims to assess the improvement provided by the branch-and-bound approach compared to what we can consider as the state-of-the-art. Furthermore, we aim to assess how the randomisation and the guided search affect the computational performances and the pruning percentage.
For the scalability analysis, we use randomly generated parameter
values and progressively increase the number of periods.  We fix the
holding cost per unit at $h=1$, but the other cost parameters are
uniform random variables: ordering cost is in the range $K \in
[80,320]$, review cost is in the range $W \in [80,320]$ and penalty
cost per unit is in the range $b \in [4,16]$.  Demands per period are
uniform random variables in the range $[30,70]$.  We generate 100
different instances and for planning horizons in the range 4--20
periods.

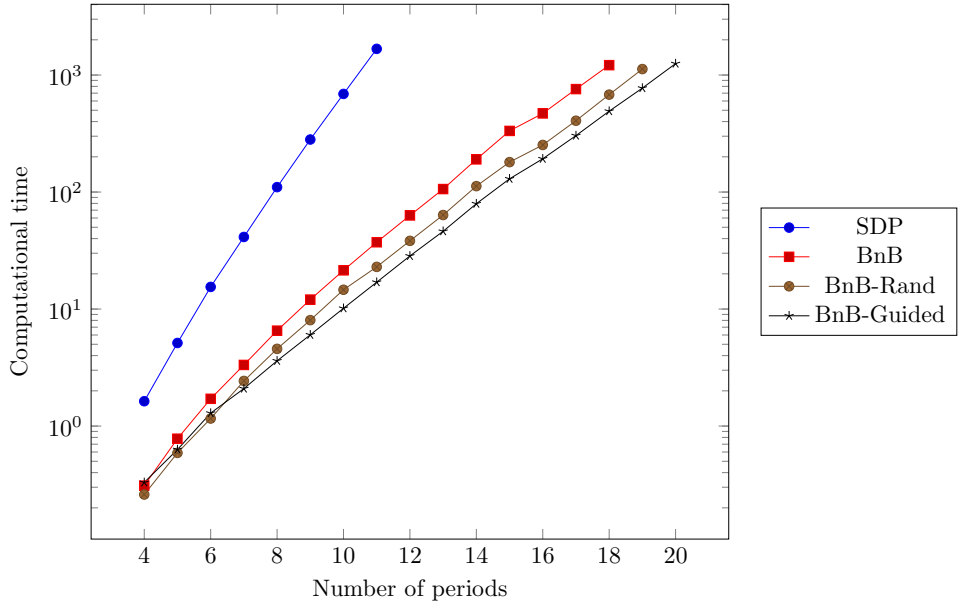
\begin{figure}
    \centering
\begin{tikzpicture} [scale = 0.9]
\pgfplotsset{every axis legend/.append style={at={(1.05,0.5)}, anchor=west}}
\begin{semilogyaxis}[
    xlabel={Number of periods},
    ylabel={Computational time}
]
\addplot coordinates{
(4,1.63)	(5,5.13)	(6,15.46)	(7,41.31)	(8,110.15)	(9,281.19)	(10,689.54)	(11,1674.18)
};
\addplot coordinates{
(4,0.31)	(5,0.78)	(6,1.71)	(7,3.33)	(8,6.54)	(9,12.05)	(10,21.42)	(11,37.23)	(12,63.07)	(13,106.04)	(14,190.2)	(15,333.01)	(16,469.6)	(17,758.17)	(18,1214.61)	
};
\addplot coordinates{
(4,0.26)	(5,0.59)	(6,1.16)	(7,2.43)	(8,4.57)	(9,8.04)	(10,14.62)	(11,22.93)	(12,38.27)	(13,63.61)	(14,112.31)	(15,180.15)	(16,252.56)	(17,406.99)	(18,680.07)	(19,1124.43)	
};
\addplot coordinates{
(4,0.33)	(5,0.63)	(6,1.29)	(7,2.09)	(8,3.61)	(9,6.02)	(10,10.15)	(11,16.96)	(12,28.43)	(13,46.2)	(14,79.44)	(15,129.65)	(16,192.41)	(17,303.55)	(18,491.09)	(19,773.62)	(20,1251.98)	
};
\legend{SDP,BnB,BnB-Rand,BnB-Guided}
\end{semilogyaxis}
\end{tikzpicture}
\caption{Average computational time of the 100 instances over the number of periods.  Time limit 1
hour.}
\label{fig:comp_time_log}
\end{figure}

Figure \ref{fig:comp_time_log} shows the average computational time over the 100 instances. The
y-axis is logarithmic to show the exponential behaviour of the
solutions.  The new method is able to solve instances almost twice as
large in a reasonable time.  Though it still has exponential
behaviour, its slope is considerably less than that of the SDP. The randomisation reduces the computational effort needed. The guided search requires the computation of an $(R,S)$ policy before the BnB approach. For small instances, the added computational effort is higher than the improvement provided by a higher pruning percentage. However, for medium/big instances, the improvement is considerable.

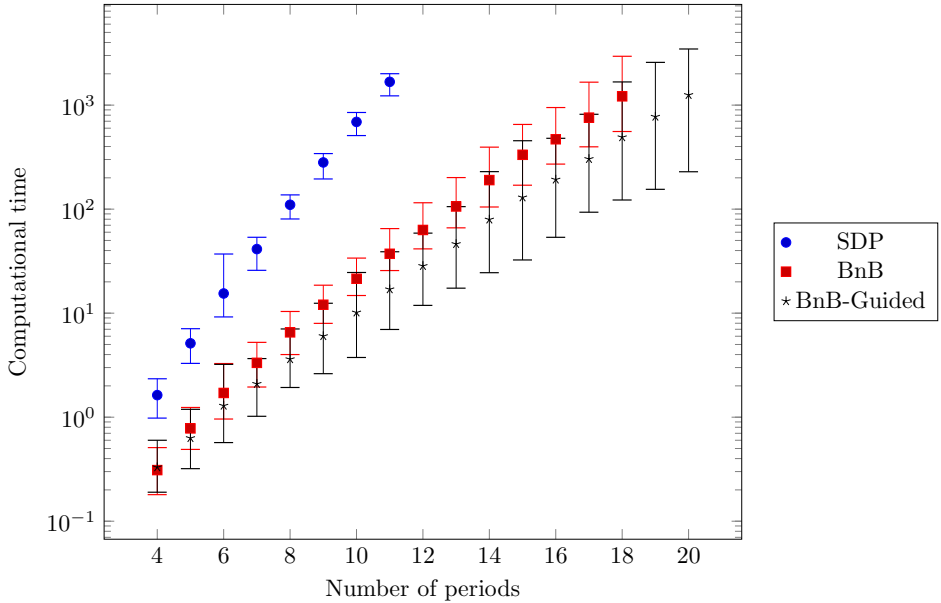
\begin{figure}
	\centering
\begin{tikzpicture}[scale = 0.9]
\pgfplotsset{every axis legend/.append style={at={(1.05,0.5)}, anchor=west}}
\begin{semilogyaxis}[xlabel={Number of periods}, ylabel={Computational time}]
\addplot+[only marks, error bars/.cd, y dir=both, y explicit,
    error bar style={line width=0pt},
    error mark options={
      rotate=90,
      blue,
      mark size=4pt,
      line width=0.5pt
    }]  coordinates{
(4,1.63) +=(0,0.71) -=(0,0.65)	(5,5.13) +=(0,1.96) -=(0,1.84)	(6,15.46) +=(0,21.58) -=(0,6.26)	(7,41.31) +=(0,12.33) -=(0,15.5)	(8,110.15) +=(0,26.93) -=(0,29.71)	(9,281.19) +=(0,60.91) -=(0,86.15)	(10,689.54) +=(0,159.88) -=(0,179.92)	(11,1674.18) +=(0,330.11) -=(0,446.69)	
};
\addplot+[only marks, error bars/.cd, y dir=both, y explicit,
    error bar style={line width=0pt},
    error mark options={
      rotate=90,
      red,
      mark size=4pt,
      line width=0.5pt
    }]  coordinates{
(4,0.31) +=(0,0.2) -=(0,0.13)	(5,0.78) +=(0,0.46) -=(0,0.29)	(6,1.71) +=(0,1.57) -=(0,0.75)	(7,3.33) +=(0,1.91) -=(0,1.38)	(8,6.54) +=(0,3.84) -=(0,2.54)	(9,12.05) +=(0,6.52) -=(0,4.07)	(10,21.42) +=(0,12.42) -=(0,6.64)	(11,37.23) +=(0,27.78) -=(0,11.55)	(12,63.07) +=(0,52.09) -=(0,21.58)	(13,106.04) +=(0,95.07) -=(0,39.94)	(14,190.2) +=(0,204.21) -=(0,85.34)	(15,333.01) +=(0,318.97) -=(0,163.15)	(16,469.6) +=(0,478.1) -=(0,198.56)	(17,758.17) +=(0,903.58) -=(0,360.8)	(18,1214.61) +=(0,1733.3) -=(0,657.07)	
};
\pgfplotsset{cycle list shift=1}

\addplot+[only marks, error bars/.cd, y dir=both, y explicit,
    error bar style={line width=0pt},
    error mark options={
      rotate=90,
      black,
      mark size=4pt,
      line width=0.5pt
    }]  coordinates{
(4,0.33) +=(0,0.27) -=(0,0.14)	(5,0.63) +=(0,0.56) -=(0,0.31)	(6,1.29) +=(0,1.93) -=(0,0.72)	(7,2.09) +=(0,1.57) -=(0,1.07)	(8,3.61) +=(0,3.44) -=(0,1.68)	(9,6.02) +=(0,6.43) -=(0,3.4)	(10,10.15) +=(0,14.45) -=(0,6.4)	(11,16.96) +=(0,21.98) -=(0,9.99)	(12,28.43) +=(0,30.37) -=(0,16.55)	(13,46.2) +=(0,59.5) -=(0,28.81)	(14,79.44) +=(0,149.9) -=(0,54.95)	(15,129.65) +=(0,324.96) -=(0,97.18)	(16,192.41) +=(0,287.28) -=(0,138.8)	(17,303.55) +=(0,512.66) -=(0,210.13)	(18,491.09) +=(0,1179.31) -=(0,368.68)	(19,773.62) +=(0,1801.87) -=(0,618.56)	(20,1251.98) +=(0,2212.14) -=(0,1023.21)	
};
\legend{SDP,BnB,BnB-Guided}
\end{semilogyaxis}
\end{tikzpicture}
\caption{Range of computational time over the number of periods.  Time limit 1 hour.}
\label{fig:comp_time_range}
\end{figure}
Figure \ref{fig:comp_time_range} shows the range of the minimum and maximum computational times for increasing planning horizon lengths; we omitted BnB-Rand to improve the readability of the plot. The SDP solution has a low variability in the required computational time. BnB-Guided presents the highest variability among the different solutions. This is due to the fact that in some instances, the pre-computed replenishment plan is the optimal one and leads to a strong pruning of the tree that reduces the computational time considerably.

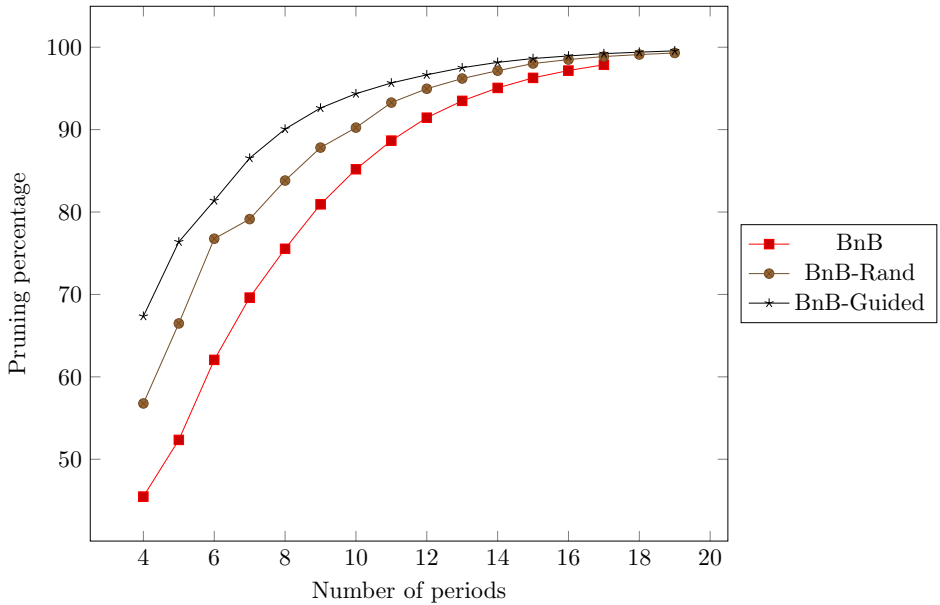
\begin{figure}
    \centering
\begin{tikzpicture}[scale = 0.9]
\pgfplotsset{every axis legend/.append style={at={(1.02,0.5)}, anchor=west}}
\begin{axis}[xlabel={Number of periods}, ylabel={Pruning percentage}]
\pgfplotsset{cycle list shift=1}

\addplot coordinates{
(4,45.48)    (5,52.35)    (6,62.06)    (7,69.62)    (8,75.54)    (9,80.93)    (10,85.19)    (11,88.67)    (12,91.45)    (13,93.49)    (14,95.06)    (15,96.28)    (16,97.17)    (17,97.87)    
};
\addplot coordinates{
(4,56.77)    (5,66.48)    (6,76.76)    (7,79.14)    (8,83.82)    (9,87.82)    (10,90.24)    (11,93.28)    (12,94.96)    (13,96.2)    (14,97.15)    (15,98.02)    (16,98.5)    (17,98.87)    (18,99.11)    (19,99.3)    
};
\addplot coordinates{
(4,67.35)    (5,76.38)    (6,81.39)    (7,86.54)    (8,90.07)    (9,92.62)    (10,94.36)    (11,95.66)    (12,96.65)    (13,97.51)    (14,98.16)    (15,98.63)    (16,98.94)    (17,99.22)    (18,99.4)    (19,99.56)    
};
\legend{BnB,BnB-Rand,BnB-Guided}
\end{axis}
\end{tikzpicture}

\caption{Average percentage of nodes pruned over the 100 instances in relation to the number of periods.}
\label{fig:prun_percentage}
\end{figure}

Figure  \ref{fig:prun_percentage}  shows the pruning percentage (Section \ref{sec:example_pruning}) of the branch-and-bound approaches.  The
pruning becomes more effective for longer planning horizons.  A high value means that the approach finds a good policy earlier in the search process; the cost of this policy provides a tighter bound for the pruning condition (Equation \ref{eq:pruning_condition_updated}). We can see that BnB-Guided provides considerable improvement.

\subsection{Instance type analysis} \label{sec:type_analysis}
In the parameter value analysis, we aim to understand how the cost parameters affect the computational effort required to find the policy and the pruning percentage. 
We use a testbed of 324 instances.
To generate the average demand values we use seasonal data with
different trends:
\begin{itemize}
    \item \textbf{(STA)} stationary case: $\tilde{d}_t = 50 $
    \item \textbf{(INC)} positive trend case: $\tilde{d}_t = \left \lceil   100t/(n-1) \right \rceil $
    \item \textbf{(DEC)} negative trend case: $\tilde{d}_t = \left \lceil  100 - 100t/(n-1) \right \rceil $
    \item \textbf{(LCY1)} life-cycle trend 1 case: this pattern is a combination of the first 3 trends. The first third of positive trend up to an average demand of 75, a central stationary one and the last negative third. If the number of periods is not a multiple of 3, the central period is extended.
    \item \textbf{(LCY2)} life-cycle trend 2 case: this pattern is a combination of INC and DEC trends. Positive trend for the first half of the planning horizon and negative trend for the second half.
    \item \textbf{(RAND)} erratic: $\tilde{d}_t = \left \lceil  U(1,100)\right \rceil $
\end{itemize}
All the patterns have an average demand of 50 per period.  For the
cost parameters we use all possible combinations of ordering cost
values $K \in \{80,160,320\}$, review costs $W \in \{ 80, 160, 320 \}$
and penalty costs $b \in \{ 4, 8, 16 \}$, with holding cost fixed at
$h = 1$.  We use all combinations of cost parameters and the six
demand patterns presented above for a full factorial experiment.  We
analyze the results for the 10-periods and 20-periods instances.

Since the baseline (equation \ref{eq:baseline} in Section
\ref{sec:sdp_ss}) is too computationally expensive (it takes approximately 45 days to solve a 20 periods instance) we replace it with an estimate in the 20-periods instances. The estimate is computed by solving 100 times the SDP for different $\gamma$ assignments and averaging it over all the possible assignments.

Tables \ref{tab:10_periods} and \ref{tab:20_periods} give an overview
of the computational time, the pruning percentage and the average
number of reviews of the methods for the 10- and 20-period
experiments.  They show that SDP is not strongly affected by the cost
parameters and that the main difference is caused by the demand
patterns.  This is due to the maximum average demand per period being
lower for STA, LCY1 and RAND.  The stationary case is faster to
compute as its maximum is 50, the second-fastest is the first life
cycle with a maximum of 75, and the erratic pattern is slowest.  All
the other patterns have a maximum of 100.


The pruning percentage gives an indication of the efficacy of the
branch-and-bound.  Our algorithms perform particularly well on high
review costs.  For instance, with 20 periods and $W = 320$ the pruning
percentage reached an impressive average of $99.83\%$ for the BnB-Guided, solving one instance in less than 13 minutes on average, while the
baseline is expected to take more than six weeks.  For the BnB,
the percentage is $98.52\%$, so it visits more than twice the nodes
compared to the guided version. The randomised search (not shown in the table for the sake of readability) reaches an average of $98.92\%$. We note that the penalty cost
also affects performance: a higher penalty cost reduces pruning.

The average number of review moments of the optimal policies decrease
as the ordering and the review increase. Also, a higher penalty cost
leads to more frequent reviews, which reduces the probability of
demand excess and mitigates the uncertainty of the inventory level. 
We observe that the decreasing pattern requires fewer review periods
than the others, due to its decreasing tail that reduces the number of
orders needed.

Our best-proposed method outperforms the baseline by factors of 50 and
1300 on 10- and 20-period instances, respectively.

\begin{table*}
\centering
\resizebox{\hsize}{!}{%
\def\arraystretch{1.5}
\renewcommand{\tabcolsep}{1.5mm}
\begin{tabular}{|c c||c|c|c||c|c||c|}
\hline
& & \multicolumn{3}{c||}{Computational time} & \multicolumn{2}{c||}{Pruning \%} &  \\ \hline
& & Base & BnB &  BnB-Guided & BnB & BnB-Guided & Nr. reviews \\ \hline
K values & 80 & 14.62	& 0.55	& 0.3	& 82.15(3.37)	& 91.51(3.94)	& 3.0\\ 
		 & 160 & 14.83	& 0.56	& 0.28	& 81.79(3.38)	& 92.06(4.07)	& 2.56\\ 
		 & 320 & 14.97	& 0.61	& 0.32	& 80.31(4.74)	& 91.06(5.37)	& 2.06\\ 
\hline
W values & 80 & 14.83	& 0.68	& 0.46	& 78.36(4.37)	& 86.94(4.08)	& 3.0\\ 
		 & 160 & 14.79	& 0.56	& 0.27	& 81.94(2.84)	& 92.48(2.48)	& 2.56\\ 
		 & 320 & 14.8	& 0.48	& 0.17	& 83.96(1.99)	& 95.2(1.76)	& 2.06\\ 
\hline
b values & 4 & 14.89	& 0.57	& 0.28	& 81.48(3.32)	& 92.06(4.21)	& 2.39\\ 
		 & 8 & 14.81	& 0.57	& 0.29	& 81.69(3.72)	& 91.89(4.27)	& 2.56\\ 
		 & 16 & 14.71	& 0.58	& 0.33	& 81.09(4.7)	& 90.68(4.93)	& 2.67\\ 
\hline
Pattern & STA & 10.76	& 0.37	& 0.21	& 82.87(2.87)	& 91.58(4.21)	& 2.63\\ 
		 & INC & 17.3	& 0.69	& 0.46	& 81.27(3.27)	& 88.67(3.96)	& 2.7\\ 
		 & DEC & 17.39	& 0.66	& 0.27	& 81.5(2.63)	& 93.7(3.6)	& 2.33\\ 
		 & LCY1 & 15.03	& 0.65	& 0.32	& 78.61(3.49)	& 90.72(4.52)	& 2.59\\ 
		 & LCY2 & 16.56	& 0.71	& 0.36	& 78.99(3.25)	& 90.57(5.01)	& 2.48\\ 
		 & RAND & 11.79	& 0.35	& 0.17	& 85.27(3.86)	& 94.01(3.19)	& 2.48\\ 
\hline
Average & & 14.81	& 0.57	& 0.3	& 81.42(3.96)	& 91.54(4.52)	& 2.54\\ 
\hline
\end{tabular}
}
\caption{Computational times (in minutes), pruning percentage and number of reviews for 10 periods instances. Between brackets is the standard deviation of the pruning percentage.}
\label{tab:10_periods}
\end{table*}
\begin{table*}
\centering
\resizebox{\hsize}{!}{%
\def\arraystretch{1.5}
\renewcommand{\tabcolsep}{1.5mm}
\begin{tabular}{|c c||c|c|c||c|c||c|}
\hline
& & \multicolumn{3}{c||}{Computational time} & \multicolumn{2}{c||}{Pruning \%} &  \\ \hline
& & Base & BnB &  BnB-Guided & BnB & BnB-Guided & Nr. reviews \\ \hline
K values & 80 & 65366.67	& 105.12	& 47.57	& 98.56(0.76)	& 99.34(0.52)	& 6.04\\ 
		 & 160 & 65470.02	& 109.35	& 49.78	& 98.53(0.92)	& 99.33(0.59)	& 5.17\\ 
		 & 320 & 66070.17	& 115.98	& 51.9	& 98.47(1.04)	& 99.32(0.68)	& 4.13\\ 
\hline
W values & 80 & 66737.03	& 181.66	& 100.37	& 97.61(0.9)	& 98.67(0.56)	& 6.04\\ 
		 & 160 & 64772.93	& 96.12	& 36.09	& 98.68(0.45)	& 99.5(0.23)	& 5.17\\ 
		 & 320 & 65396.9	& 52.67	& 12.8	& 99.27(0.23)	& 99.83(0.1)	& 4.13\\ 
\hline
b values & 4 & 65851.88	& 96.59	& 41.92	& 98.7(0.76)	& 99.45(0.53)	& 4.78\\ 
		 & 8 & 65847.45	& 108.37	& 48.75	& 98.56(0.83)	& 99.35(0.55)	& 5.2\\ 
		 & 16 & 65207.52	& 125.49	& 58.59	& 98.3(1.07)	& 99.21(0.69)	& 5.35\\ 
\hline
Pattern & STA & 43447.11	& 73.24	& 38.54	& 98.51(0.84)	& 99.23(0.66)	& 5.3\\ 
		 & INC & 72449.66	& 110.73	& 69.98	& 98.69(0.62)	& 99.18(0.53)	& 5.41\\ 
		 & DEC & 72706.98	& 141.49	& 48.06	& 98.29(1.05)	& 99.43(0.57)	& 4.7\\ 
		 & LCY1 & 62607.87	& 139.2	& 62.58	& 98.05(0.96)	& 99.13(0.72)	& 5.19\\ 
		 & LCY2 & 69243.25	& 141.35	& 56.22	& 98.22(0.82)	& 99.3(0.57)	& 5.04\\ 
		 & RAND & 73358.85	& 54.88	& 23.13	& 99.36(0.31)	& 99.74(0.23)	& 5.04\\ 
\hline
Average & & 65635.62	& 110.15	& 49.75	& 98.52(0.91)	& 99.33(0.6)	& 5.11\\ 
\hline
\end{tabular}
}
\caption{Computational times (in minutes), pruning percentage and number of reviews for 20 periods instances. Between brackets the standard deviation of the pruning percentage.}
\label{tab:20_periods}
\end{table*}

\section{Conclusion and Future Work} \label{sec:conclusion}

In this paper, we considered a single-item single-stocking location
inventory lot-sizing problem with non-stationary stochastic demand,
fixed and linear ordering cost, review cost, holding cost and penalty cost.  We
present the first algorithm to compute optimal $(R,s,S)$ policy
parameters.  This policy has a high practical value, but the
computation of optimal or near-optimal parameters has been considered
extremely difficult.  Our proposed technique is a hybrid of
branch-and-bound and stochastic dynamic programming, enhanced by ad
hoc bounds computed with dynamic programming, and by a randomised
depth-first exploration of the search tree.

In an extensive numerical study, we first investigated the scalability
of the technique at increasing time horizon, analysing both
computational time and the efficacy of the bounding technique.  We
then tested the performance of the method for different cost
parameters.  Our technique performs best on low penalty costs and high
review costs.  On 20-period instances, it outperforms a baseline method
by three orders of magnitude.

This technique opens up multiple research directions on the
determination of $(R,s,S)$ policy parameters.  It can lead to new
optimal solutions for the same problem, and it can be improved with
tighter bounds.  It is also useful for computing optimality gaps of
new heuristics.

As future studies, the approach presented herein can be extended to overcome some of the limitations of the problem setting. Considering multiple items with joint shipping or modelling a more complex supply chain with multiple echelons are generalisations of this problem, and they would increase the applicability of the $(R,s,S)$ policy.

\subsubsection*{Acknowledgments}
This publication has emanated from research conducted with the financial support of Science Foundation Ireland under Grant number 12/RC/2289-P2, which is co-funded under the European Regional Development Fund.

\bibliography{mybibfile}

\end{document}